\documentclass[11pt]{amsart}
\usepackage{amscd,amssymb,amsmath,hyperref}
\input xy
\xyoption{all}

\newtheorem{theorem}{Theorem}[section]

\newtheorem{conjecture}[theorem]{Conjecture}
\newtheorem{question}[theorem]{Question}

\renewcommand{\leq}{\leqslant}
\renewcommand{\geq}{\geqslant}

\theoremstyle{definition}
\newtheorem{definition}[theorem]{Definition}
\newtheorem{example}[theorem]{Example}

\theoremstyle{definition}

\numberwithin{equation}{section}

\newcommand{\ve}{\varepsilon}
\newcommand{\vp}{\varphi}

\newcommand{\ov}[1]{\overline{#1}}
\newcommand{\de}{\partial}
\newcommand{\db}{\overline{\partial}}
\newcommand{\ddbar}{\sqrt{-1} \partial \overline{\partial}}
\newcommand{\ti}[1]{\tilde{#1}}
\newcommand{\dVol}{\mathrm{dVol}}

\newcommand{\tr}[2]{\textrm{tr}_{#1} #2}
\DeclareMathOperator{\Ric}{Ric}

\numberwithin{equation}{section} \numberwithin{figure}{section}

\author{Valentino Tosatti}
\address{Department of Mathematics, Northwestern University, 2033 Sheridan Road, Evanston, IL 60208}
\email{tosatti@math.northwestern.edu}
\title{Ricci-flat metrics and dynamics on K3 surfaces}

\begin{document}

\begin{abstract}We give an overview of some recent interactions between the geometry of $K3$ surfaces and their Ricci-flat K\"ahler metrics and the dynamical study of $K3$ automorphisms with positive entropy.
\end{abstract}

\maketitle

\section{Introduction}
$K3$ surfaces form a distinguished class of compact complex surfaces which has received a tremendous amount of attention in several branches of mathematics. Our interest in $K3$ surfaces stems from the fact that they are $2$-dimensional Calabi-Yau manifold and hence admit Ricci-flat (but not flat) K\"ahler metrics, as we will explain below. The geometry of these metrics is still not completely understood, especially when families of such metrics degenerate. In a seemingly unrelated direction, $K3$ surfaces have also been studied in holomorphic dynamics. The theory of holomorphic dynamics in $1$ complex variable (on the Riemann sphere) is of course an enormous research area, and when one passes to $2$ complex variables, it turns out that the only dynamically interesting automorphisms exist on $K3$ and rational surfaces (see \cite{Can2} for the precise statement), and interesting $K3$ automorphism are relatively easy to construct. The dynamical study of such automorphisms was initiated by Cantat \cite{Can}, and we refer the reader to the survey articles \cite{CanS,CanS2,CanS3} and lecture notes \cite{Fil} for a broader overview.

The goal of this article will be to give an introduction to both of these aspects related to $K3$ surfaces and to explain some recent work by Filip and the author \cite{FT,FT2} that exploits Ricci-flat metrics to prove results in dynamics and vice versa.

In section \ref{si} we give an introduction to $K3$ surfaces, including basic examples, the conjectures of Andreotti and Weil and their solutions. In section \ref{sy} we discuss Yau's Theorem \cite{Ya} on the existence of Ricci-flat K\"ahler metrics on $K3$ surfaces. Section \ref{sd} gives an overview of the dynamical study of automorphisms of $K3$ surfaces, including basic properties and examples. Section \ref{s1} discusses the recent Kummer rigidity theorem of Cantat-Dupont \cite{CD} and Filip and the author \cite{FT2}, with an emphasis on the application of Ricci-flat metrics to this result that was found in \cite{FT2}. In section \ref{s2} we discuss applications of dynamics (in particular of Kummer rigidity) to the problem of understanding the behavior of Ricci-flat K\"ahler metrics on $K3$ surfaces when the K\"ahler class degenerates, following our work in \cite{FT}. Lastly, in section \ref{sc} we discuss a few related open problems.\\

{\bf Acknowledgements. }This article is based on the author's plenary lecture at the XXI Congress of the Italian Mathematical Union that took place in Pavia on September 2--7, 2019, and the author would like to express his gratitude to the Scientific Committee for the invitation and to all the Organizers of this wonderful Congress. He is also grateful to S. Filip for many useful discussions and for helpful comments on this paper. The author was partially supported by NSF grant DMS-1903147. This article was written during the author's stay at the Center for Mathematical Sciences and Applications at Harvard University, which he would like to thank for the hospitality.

\section{$K3$ Surfaces}\label{si}
\subsection{Complex manifolds}
The main object of study in this article are $K3$ surfaces (over the complex numbers). Before we get to the definition, let us briefly recall the basic definition of complex manifold. A complex $n$-dimensional manifold is a real manifold $X$ (of real dimension $2n$) which admits an atlas with charts with values in $\mathbb{C}^n\cong\mathbb{R}^{2n}$ whose transition maps are {\em holomorphic} (i.e. complex analytic). We will implicitly assume that all our manifolds are Hausdorff, second-countable and connected.

The first examples of complex manifolds are Riemann surfaces, which are $1$-dimensional complex manifolds. Some basic examples of compact complex manifolds include complex tori $X=\mathbb{C}^n/\Lambda,$ where $\Lambda\cong\mathbb{Z}^{2n}$ is a lattice in $\mathbb{C}^n$, complex projective space $\mathbb{CP}^n=(\mathbb{C}^{n+1}\backslash\{0\})/\mathbb{C}^*$ (acting diagonally), and smooth projective varieties $X=\{P_1=\cdots=P_m=0\}\subset\mathbb{CP}^N$, where the $P_j$'s are homogeneous polynomials (and we assume of course that $X$ is a manifold). On the other hand, a compact complex manifold is called projective if it admits a holomorphic embedding into $\mathbb{CP}^N$ for some $N$, and thanks to a classical theorem of Chow its image is cut out by finitely many polynomial equations, thus giving us a smooth projective variety.

\subsection{$K3$ surfaces}Complex tori $X=\mathbb{C}^n/\Lambda$ have a special property: they admit a never-vanishing holomorphic $n$-form  $\Omega$, which is induced by $dz_1\wedge\cdots\wedge dz_n$ on $\mathbb{C}^n$ (which is obviously translation-invariant). It is a classical result in the theory of Riemann surfaces that when $n=1$ the property of admitting a never-vanishing holomorphic $1$-form characterizes $1$-dimensional tori (elliptic curves) among compact Riemann surfaces. However, this is not true anymore for compact complex manifolds of dimension $n\geq 2$, and indeed our main object of study is the following:

\begin{definition}
A compact complex manifold $X$ of complex dimension $2$ is called a $K3$ surface if $X$ is simply connected and it admits a never-vanishing holomorphic $2$-form $\Omega$.
\end{definition}

The first examples of $K3$ surfaces were studied in the $19$th century by Kummer, Cayler, Schur and others, and later by the Italian school of algebraic geometry, in particular by Enriques and Severi. After work of Andreotti and Atiyah in the early 1950s, these surfaces were given their name by Weil \cite{We}, who in a grant report laid out four basic conjectures about $K3$ surfaces, that were also independently formulated by Andreotti, and that shaped the research in the field for the coming decades.

Before we get to that, let us give some basic examples of $K3$ surfaces. The reader is refereed to the classic textbooks \cite{BHPV,K3} and the more recent \cite{Huy,Kon} for details.
\subsection{Examples}
\begin{example}[Quartic surfaces] Every smooth hypersurface $X=\{P=0\}\subset\mathbb{CP}^3$ with $\deg P=4$ is a $K3$ surface. Perhaps the simplest quartic surface is the Fermat quartic given by
$$z_0^4+z_1^4+z_2^4+z_3^4=0.$$
\end{example}

\begin{example}[Complete intersections in products of projective spaces]
More generally, we can consider smooth complete intersections in a product of $k$ projective spaces,
$$X=\{P_1=\cdots=P_m=0\}\subset\mathbb{CP}^{n_1}\times\cdots\times \mathbb{CP}^{n_k},$$
where each $P_j$ is a multihomogeneous polynomial (i.e. homogeneous separately in the homogeneous coordinates of each of the $\mathbb{CP}^{n_p}$ factors) of multidegree $\deg P_j=(d_1^{(j)},\dots,d_k^{(j)})$ so that we have $\sum_{p=1}^k n_p=m+2$ (hence $X$ is complex $2$-dimensional) and
$$\sum_{j=1}^m d_p^{(j)}=n_p+1,$$
for all $1\leq p\leq k$.

Some explicit examples are:
\begin{itemize}
\item The complete intersection of a quadric and a cubic in $\mathbb{CP}^4$ ($k=1, m=2, d_1^{(1)}=2, d_1^{(2)}=3$)
\item The complete intersection of three quadrics in $\mathbb{CP}^5$ ($k=1,m=3, d_1^{(1)}=d_1^{(2)}=d_1^{(3)}=2$)
\item Smooth hypersurfaces in $\mathbb{CP}^2\times\mathbb{CP}^1$ of multidegree $(3,2)$ ($k=2,n_1=2,n_2=1,m=1, d_1^{(1)}=3,d_2^{(1)}=2$)
\item Smooth hypersurfaces in $\mathbb{CP}^1\times\mathbb{CP}^1\times\mathbb{CP}^1$ of multidegree $(2,2,2$) ($k=3,n_1=n_2=n_3=1,m=1,d_1^{(1)}=d_2^{(1)}=d_3^{(1)}=2$)
\item Complete intersections of two hypersurfaces of bidegrees $(1,1)$ and $(2,2)$ respectively in $\mathbb{CP}^2\times\mathbb{CP}^2$ ($k=2, n_1=n_2=2,m=2, d_1^{(1)}=d_2^{(1)}=1,d_1^{(2)}=d_2^{(2)}=2$).
\end{itemize}
\end{example}

\begin{example}[Kummer surfaces]
Here we start with a $2$-torus $T=\mathbb{C}^2/\Lambda$ and consider the involution $\iota$ of $T$ which is induced by $\iota(z_1,z_2)=(-z_1,-z_2)$. The involution has $16$ fixed points, which become $16$ singularities of the quotient $Y=T/\iota$. These singularities are rational double points (orbifold points with orbifold group $\mathbb{Z}/2\mathbb{Z}$) which can be resolved by a simple blowup, to obtain $\pi:X\to Y$ where $X$ is a smooth compact complex surface. It is not hard to verify (see e.g. \cite{BHPV,K3,Huy,Kon}) that $X$ is $K3$, the Kummer surface associated to $T$. The preimage under $\pi$ of the $16$ singular points of $Y$ are $16$ rational curves in $X$ with self-intersection $-2$.
\end{example}
\subsection{The conjectures of Andreotti and Weil}As mentioned above, in the 1950s the study of $K3$ surfaces shifted its focus from specific examples to a general theory. The following $4$ basic conjectures were made independently by Andreotti and Weil \cite{We}:\\

(I) All $K3$ surfaces form one connected (complex-analytic) family. In particular they are all diffeomorphic the same smooth $4$-manifold. This conjecture was proved by Kodaira \cite{Ko} in 1964.\\

It is interesting to remark that the family of projective $K3$ surfaces is $19$-dimensional (this is the same dimension as the space of smooth quartics in $\mathbb{CP}^3$) while the family of all $K3$ surfaces is $20$-dimensional, so in a sense most $K3$ surfaces are not projective.

For every $K3$ surface $X$, the cohomology $H^2(X,\mathbb{Z})$ equipped with the intersection form is isomorphic to a fixed lattice $\Lambda$, which is the unique even unimodular lattice of signature $(3,19)$. A marking on $X$ is then a choice of isomorphism of lattices $\iota:\Lambda\to H^2(X,\mathbb{Z})$.
Now on $X$ the never-vanishing holomorphic $2$-form $\Omega$ is unique up to scaling, and it satisfies $\Omega\wedge\Omega=0$ while $\Omega\wedge\ov{\Omega}$ is a smooth positive volume form on $X$, so that $\int_X\Omega\wedge\ov{\Omega}>0$. Thus, if we denote also by $\iota:\Lambda\otimes\mathbb{C}\to H^2(X,\mathbb{C})$ the induced isomorphism, then $\iota^{-1}([\Omega])$ gives a well-defined point $\mathcal{P}(X,\iota)$ in the period domain
$$\mathcal{D}=\{c\in\mathbb{P}(\Lambda\otimes\mathbb{C})\ |\ c\cdot c=0,\quad c\cdot \ov{c}>0\}.$$
The map $\mathcal{P}$ that associates to a marked $K3$ surface $(X,\iota)$ its period point $\mathcal{P}(X,\iota)$ is called the period map. The second conjecture is then:\\

(II) (Torelli Theorem) If two marked $K3$ surfaces $(X,\iota),(X',\iota')$ determine the same period point $\mathcal{P}(X,\iota)=\mathcal{P}(X',\iota')$, then $X$ and $X'$ are biholomorphic. This conjecture was proved by Pjatecki\u{i}-\v{S}apiro-\v{S}afarevi\v{c} \cite{PS} in 1971 for projective $K3$ surfaces and by Burns-Rapoport \cite{BR} in 1975 in general.\\

To state the next conjecture, we need another basic definition. A Hermitian metric $g$ on a complex manifold $X^n$ is a smoothly-varying family of Hermitian inner products on the holomorphic tangent spaces $T^{1,0}_x X$, which in local holomorphic coordinates is thus given by an $n\times n$ positive definite Hermitian matrix $(g_{j\ov{k}}(x))_{j,k=1}^n$ which varies smoothly in $x$. Every complex manifold admits Hermitian metrics, as can be seen for example by patching together local Euclidean inner products using a partition of unity. To a Hermitian metric $g$ one then associates a smooth real $(1,1)$-form $\omega$ which in local coordinates is given by
$$\omega=i\sum_{j,k=1}^ng_{j\ov{k}}dz_j\wedge d\ov{z}_k,$$
and we say that $g$ (or $\omega$) is K\"ahler if $d\omega=0$. This can be viewed as a system of first-order linear PDEs for the coefficients $g_{j\ov{k}}$, and the existence of a K\"ahler metric on a compact complex manifold implies several nontrivial global constraints (for example the even Betti numbers of $X$ must be nonzero, and the odd Betti numbers must be even). The restriction of a K\"ahler metric to a complex submanifold is still K\"ahler, and since $\mathbb{CP}^n$ admits the explicit Fubini-Study K\"ahler metric, it follows that all projective manifolds admit K\"ahler metrics.\\

(III) Every $K3$ surface admits a K\"ahler metric. This was proved by Siu \cite{Si} in 1983.\\

Combined with earlier work of Kodaira \cite{Ko} and Miyaoka \cite{Mi}, this result implies that a compact complex surface is K\"ahler if and only if its first Betti number is even. New proofs of this result were later found by Buchdahl \cite{Bu} and Lamari \cite{La} independently.\\

(IV) (Surjectivity of the period map) Every point in $\mathcal{D}$ is the period point of some marked $K3$ surface. This was proved by Kulikov \cite{Ku} in 1977 for projective $K3$ surfaces and by Todorov \cite{Tod} in 1980 in general.\\

Lastly, we mention that there is a refined Torelli Theorem, which is also proved in \cite{BR,PS}. For a $K3$ surface $X$, the set $\mathcal{C}_X$ of all cohomology classes of K\"ahler metrics is a cone inside $H^2(X,\mathbb{R})$. The refined Torelli Theorem then states that if two marked $K3$ surfaces $(X,\iota),(X',\iota')$ determine the same period point $\mathcal{P}(X,\iota)=\mathcal{P}(X',\iota')$ and furthermore $\iota\circ\iota'^{-1}$ takes $\mathcal{C}_{X'}$ to $\mathcal{C}_X$, then there is a unique biholomorphism $F:X\to X'$ such that $F^*=\iota\circ\iota'^{-1}$.

\section{Ricci-flat K\"ahler metrics}\label{sy}
\subsection{Ricci curvature}Thanks to the aforementioned theorem of Siu, every $K3$ surface admits a K\"ahler metric. As we will now see, they admit rather special K\"ahler metrics.

Recall that given a K\"ahler metric $g$, one has the associated Ricci curvature tensor $R_{j\ov{k}}$, which as in Riemannian geometry is the trace of the Riemann curvature tensor. The fact that $g$ is K\"ahler implies that the Ricci tensor is Hermitian (i.e. $\ov{R_{j\ov{k}}}=R_{k\ov{j}}$) and that the associated real $(1,1)$-form
$$\Ric(g)=i\sum_{j,k=1}^nR_{j\ov{k}}dz_j\wedge d\ov{z}_k,$$
is closed, $d\Ric(g)=0$, and it is locally given by
$$\Ric(g)=-i\de\db\log\det(g_{j\ov{k}}),$$
where $d=\de+\db$ is the usual splitting of the exterior derivative on a complex manifold.

\subsection{Ricci-flatness}Let $X$ be a $K3$ surface. Since the never-vanishing holomorphic $2$-form $\Omega$ is unique up to scaling, we will assume from now on that it has been scaled so that the smooth positive volume form $\dVol:=\Omega\wedge\ov{\Omega}$ satisfies $\int_X\dVol=1$.

If $\omega$ is a K\"ahler metric on $X$, then its volume form $\omega^2$ can be written as $\omega^2=f \dVol$ for some smooth positive function $f$ on $X$. In local holomorphic coordinates we can write $\omega=ig_{j\ov{k}}dz_j\wedge d\ov{z}_k$ (using now the Einstein summation convention) and $\Omega=h dz_1\wedge dz_2$, where $h$ is a locally-defined never-vanishing holomorphic function. It then follows that
$$\det(g_{j\ov{k}})=f|h|^2,$$
and taking $-i\de\db\log$ of this, and using that $i\de\db\log|h|^2=0$ (an elementary computation), we get
$$\Ric(g)=-i\de\db\log f.$$
From this we see that if $f$ is constant then $\Ric(g)$ vanishes identically, and the converse is also true since if $i\de\db\log f=0$ then the strong maximum principle implies that $f$ is constant. Thus, $\omega$ is a K\"ahler metric with vanishing Ricci curvature if and only if its volume form is a constant multiple of $\dVol$,
\begin{equation}\label{rf}
\omega^2=c\,\dVol,\quad c\in\mathbb{R}_{>0},
\end{equation}
and of course integrating this identity we see that $c=\int_X\omega^2$.

\subsection{Yau's Theorem}Now, if we start with a K\"ahler metric $\omega$ on a $K3$ surface $X$, and write $\omega^2=f\dVol$ as above (with $f$ not necessarily constant), it is then clear that the conformally rescaled Hermitian metric $\ti{\omega}=e^{-\frac{f}{2}}\omega$ satifies $\ti{\omega}^2=\dVol$, which is the equation we want, but it is easy to see that $\ti{\omega}$ will not be closed (and so the corresponding Hermitian metric will not be K\"ahler) unless $f$ is a constant.

On the other hand, if $\ti{\omega}$ and $\omega$ are two K\"ahler metrics on $X$, they define cohomology classes $[\ti{\omega}]$ and $[\omega]$ in $H^{1,1}(X,\mathbb{R})$, the subspace of $H^2(X,\mathbb{R})$ of de Rham cohomology classes which admit a representative which is a closed real $(1,1)$-form. The basic $\de\db$-Lemma of Kodaira shows that if $[\ti{\omega}]=[\omega]$ in $H^{1,1}(X,\mathbb{R})$, then there exists $\vp\in C^\infty(X,\mathbb{R})$, unique up to an additive constant, such that $$\ti{\omega}=\omega+\ddbar\vp,$$
which in local coordinates translates to
$$\ti{g}_{j\ov{k}}=g_{j\ov{k}}+\frac{\de^2\vp}{\de z_j\de\ov{z}_k}.$$
Thus, to find a Ricci-flat K\"ahler metric $\ti{\omega}$ with $[\ti{\omega}]=[\omega]$, it suffices to find $\vp\in C^\infty(X,\mathbb{R})$ such that
\begin{equation}\label{eqn}
\omega+\ddbar\vp>0,\quad (\omega+\ddbar\vp)^2=\left(\int_X\omega^2\right)\dVol,
\end{equation}
where we used Stokes's Theorem $\int_X\ti{\omega}^2=\int_X(\omega+\ddbar\vp)^2=\int_X\omega^2$.
In local coordinates, \eqref{eqn} becomes a fully nonlinear PDE of complex Monge-Amp\`ere type
$$\left(g_{j\ov{k}}+\frac{\de^2\vp}{\de z_j\de\ov{z}_k}\right)>0\quad\det\left(g_{j\ov{k}}+\frac{\de^2\vp}{\de z_j\de\ov{z}_k}\right)=\left(\int_X\omega^2\right)|h|^2,$$
where $h$ is a local never-vanishing holomorphic function as above. The fundamental result is then:
\begin{theorem}[Yau 1976 \cite{Ya}]\label{yau}
Let $X$ be a $K3$ surface equipped with a K\"ahler metric $\omega$, and let $\Omega$ be the never-vanishing holomorphic $2$-form on $X$ normalized so that $\dVol=\Omega\wedge\ov{\Omega}$ has integral $1$. Then there exists $\vp\in C^\infty(X,\mathbb{R})$, unique up to an additive constant, such that \eqref{eqn} holds. The K\"ahler metric $\ti{\omega}=\omega+\ddbar\vp$ is then a Ricci-flat K\"ahler metric on $X$, the unique such metric with $[\ti{\omega}]=[\omega]$.
\end{theorem}
This is a special case of Yau's solution \cite{Ya} of the Calabi Conjecture \cite{Cal}, which solves an equation analogous to \eqref{eqn} on arbitrary compact K\"ahler manifolds.

\subsection{The Hyperk\"ahler property}The Ricci-flat K\"ahler metrics $\ti{g}$ that are produced by Theorem \ref{yau} are not explicit, as is often the case for solutions of nonlinear PDEs. On the other hand, it is not hard to see (see e.g. \cite{Be}) that these metrics have resticted holonomy: the holonomy group $\mathrm{Hol}(\ti{g})$ of linear transformations of $T_xX$ obtained by $\ti{g}$-parallel transport along loops based at $x$ (an arbitrary basepoint) is precisely equal to $SU(2)=Sp(1)$. This implies that the metric $\ti{g}$ are hyperk\"ahler: the manifold $X$ admits a triple of complex structures $I,J,K,$ which satisfy the quaternionic relations ($I\circ J=K$, etc.) such that $\ti{g}$ is K\"ahler with respect to $I,J$ and $K$. This last condition means that $\ti{g}$ satisfies the Hermitian property $\ti{g}(I\cdot,I\cdot)=\ti{g}(\cdot,\cdot)$ (and the same for $J,K$), and the K\"ahler form $\ti{\omega}_I(\cdot,\cdot)=\ti{g}(I\cdot,\cdot)$ is closed (and the same for $J,K$). We may assume that $I$ is the same complex structure as the one that we had fixed earlier on $X$, so that with our notation we have $\ti{\omega}=\ti{\omega}_I$, and then if $\Omega$ is as before a holomorphic $2$-form on $X$ (with the complex structure $I$), then after suitable rescaling we have $\ti{\omega}_J=\mathrm{Re}\,\Omega$ and $\ti{\omega}_K=\mathrm{Im}\,\Omega$.

It follows immediately that $\ti{g}$ is also K\"ahler with respect to all the complex structures of the form $aI+bJ+cK$ with $a,b,c\in\mathbb{R},a^2+b^2+c^2=1$, which is an $S^2$-worth of complex structures called the twistor sphere of $X$. Passing from one of these complex structures to another one is usually referred to as hyperk\"ahler rotation, and while this changes the complex structure and the K\"ahler form, it does not change the Hermitian metric $\ti{g}$.

\section{Dynamics of $K3$ automorphisms}\label{sd}
\subsection{$K3$ automorphisms} We now shift our attention to the group $\mathrm{Aut}(X)$ of automorphisms (i.e. biholomorphisms) of a $K3$ surface $X$. This is in general a discrete group, in fact it embeds as a subgroup of the orthogonal group of $H^2(X,\mathbb{Z})$ equipped with the intersection form, but it can still be quite large, as we will see later (cf. Examples \ref{we} and \ref{222}).

First, let us make the following observation. Let $T:X\to X$ be an automorphism of a $K3$ surface $X$, equipped with its normalized holomorphic $2$-form $\Omega$. Then the pullback $T^*\Omega$ is also a never-vanishing holomorphic $2$-form on $X$, and so it must be a constant multiple of $\Omega$, $T^*\Omega=c\Omega, c\in\mathbb{C}$. But since $\dVol=\Omega\wedge\ov{\Omega}$ is a positive volume form, and $T$ is an automorphism, we see that
$$\int_X\Omega\wedge\ov{\Omega}=\int_XT^*(\Omega\wedge\ov{\Omega})=|c|^2\int_X\Omega\wedge\ov{\Omega},$$
so $|c|^2=1$, and therefore we see that
$$T^*\dVol=\dVol,$$
i.e. the volume form $\dVol$ is $\mathrm{Aut}(X)$-invariant.

\subsection{Hyperbolic geometry}Since $T$ is holomorphic, pullback by $T$ gives a linear map $T^*:H^{1,1}(X,\mathbb{R})\to H^{1,1}(X,\mathbb{R})$, which preserves the intersection pairing on $H^{1,1}(X,\mathbb{R})$ (which has signature $(1,19)$). Let us consider the $2$-sheeted hyperboloid
$$\{c\in H^{1,1}(X,\mathbb{R})\ |\ c\cdot c=1\}.$$
Then $T^*$ preserves this hyperboloid, and it also preserves its two sheets. Let $\mathcal{H}$ be the sheet which contains the cohomology class $\frac{[\omega]}{\int_X\omega^2}$ of some (and hence any) K\"ahler metric $\omega$, then $\mathcal{H}$ with its intersection form is a model of hyperbolic space $\mathbb{H}^{19}$ and $T^*:\mathcal{H}\to\mathcal{H}$ gives an isometry of hyperbolic space.

Isometries of hyperbolic space can be divided into three classes: elliptic if they admit a fixed point in $\mathcal{H}$, parabolic if they admit a unique fixed point on the ideal boundary $\de\mathcal{H}$, and hyperbolic if they admit two fixed points on $\de\mathcal{H}$. We then have the following remarkable result:
\begin{theorem}[Cantat 1999 \cite{Can}]
Let $T:X\to X$ be a $K3$ automorphism, and $T^*:\mathcal{H}\to\mathcal{H}$ the corresponding isometry of hyperbolic $19$-space. Then
\begin{itemize}
\item $T^*$ is elliptic $\Leftrightarrow$ $T$ is of finite order (i.e. $T^k=\mathrm{Id}$ for some $k\geq 1$)
\item $T^*$ is parabolic $\Leftrightarrow$ $T$ is of infinite order and it preserves an elliptic fibration $\pi:X\to\mathbb{CP}^1$
\end{itemize}
\end{theorem}
Here an elliptic fibration on a $K3$ surface $X$ is a surjective holomorphic map $\pi:X\to\mathbb{CP}^1$ with connected fibers and with generic fibers elliptic curves. Such elliptic fibrations have a nonzero and finite number of singular/multiple fibers. An automorphism $T$ of $X$ is said to preserve an elliptic fibration $\pi$ if it maps every fiber of $\pi$ to itself.

\subsection{Hyperbolic automorphisms and entropy}It is then natural to ask what happens when $T^*$ is hyperbolic. From the definition it follows that $T^*$ is hyperbolic if and only if the spectral radius $\rho$ of $T^*:H^{1,1}(X,\mathbb{R})\to H^{1,1}(X,\mathbb{R})$ is strictly larger than $1$.

This spectral radius turns out to be related to a basic quantity in the study of the dynamical behavior of iterates $T^n$ of $T$, $n\geq 1$: the topological entropy. This can be defined as follows. Fix any K\"ahler metric on $X$, denote by $d$ its induced distance function on $X$, and for $n\geq 1 $ and $\ve>0$ we say that a subset $A\subset X$ is $(n,\ve)$-separated if for all distinct $x,y\in A$ there is some $0\leq j\leq n$ such that $d(T^j(x),T^j(y))>\ve$. Since $X$ is compact, every $(n,\ve)$-separated subset must be finite and we let $r(n,\ve)$ to be the maximal cardinality of an $(n,\ve)$-separated subset of $X$. We then define the topological entropy $h(T)$ of $T$ by
$$h(T)=\lim_{\ve\downarrow 0}\limsup_{n\to+\infty}\frac{1}{n}\log r(n,\ve)\in [0,+\infty).$$
It is clear that this does not depend on the choice of metric on $X$. Informally, $h(T)$ measures the exponential growth rate of distinguishable orbits of $T$ when we observe the dynamics only up to $n$ iterates. The quantity inside the $\lim_{\ve\downarrow 0}$ is the same growth rate when we are only able to make measurements with precision $\ve.$ If the entropy is strictly positive, there is a very strong ``dependence on the initial conditions'', and we will think of this as one of the incarnation of a chaotic dynamical system.

The fundamental result is then the following:
\begin{theorem}[Gromov 1976 \cite{Gr}, Yomdin 1987 \cite{Yo}]
Let $T:X\to X$ be an automorphism of a $K3$ surface, and let $\rho$ be the spectral radius of $T^*:H^{1,1}(X,\mathbb{R})\to H^{1,1}(X,\mathbb{R})$. Then the topological entropy $h(T)$ of $T$ equals
$$h(T)=\log\rho.$$
\end{theorem}
This theorem combines the inequality $h(T)\geq \log\rho$ due to Yomdin \cite{Yo} in much greater generality, and the reverse inequality by Gromov \cite{Gr}. The theorem remains true for holomorphic self-maps $T:X\to X$ of any compact K\"ahler manifold $X^n$, with $\rho$ now being the maximum for $1\leq k\leq n$ of the spectral radius of $T^*$ on $H^{k,k}(X,\mathbb{R})$.

\subsection{Examples}By the Gromov-Yomdin Theorem, hyperbolic automorphisms of $K3$ surfaces are exactly those with positive topological entropy. Here we give some examples of such automorphisms.
\begin{example}[Kummer examples]\label{ku}
Let $Y=\mathbb{C}^2/\Lambda$ be a $2$-torus with an automorphism $T_Y$ with positive topological entropy $h(T_Y)>0$. Note that the automorphism $T_Y$ is induced by an affine linear map on $\mathbb{C}^2$. Then if $X$ is the Kummer $K3$ surface associated to $Y$, then $T_Y$ lifts to an automorphism $T$ of $X$, which satisfies $h(T)=h(T_Y)$. We will refer to such $(X,T)$ as Kummer examples.

For an explicit construction, one can take $\Lambda=(\mathbb{Z}\oplus i\mathbb{Z})^2$ the ``square'' lattice, and $T_Y$ induced by ``Arnol'd's cat map'' $\begin{pmatrix}
2&1\\
1&1\end{pmatrix}.$ Then one can compute that the spectral radius equals $\rho=\left(\frac{3+\sqrt{5}}{2}\right)^2$, and so $h(T)\sim 1.92$.
\end{example}

\begin{example}[Wehler \cite{Weh}]\label{we} Consider a $K3$ surface $X$ which is a complete intersection of two general hypersurfaces of bidegrees $(1,1)$ and $(2,2)$ in $\mathbb{CP}^2\times\mathbb{CP}^2$. The two projection maps to $\mathbb{CP}^2$ exhibit $X$ as a ramified double cover of $\mathbb{CP}^2$, let $\sigma_1,\sigma_2$ be the two covering involutions, and $T=\sigma_1\circ\sigma_2$. Then $T$ is an automorphism of $X$ with $h(T)=\log\left(\frac{13+\sqrt{165}}{2}\right)\sim 2.56$ (see \cite{CaT}). It is shown in \cite{Weh} that $\mathrm{Aut}(X)\cong \mathbb{Z}/2\mathbb{Z} * \mathbb{Z}/2\mathbb{Z}$, the free product generated by these two involutions.
\end{example}

\begin{example}[Mazur \cite{Maz}]\label{222} Let now $X$ be a $K3$ surface with is a generic hypersurface of multidegree $(2,2,2)$ in $\mathbb{CP}^1\times\mathbb{CP}^1\times\mathbb{CP}^1$. Now we have three projection maps to $\mathbb{CP}^1\times\mathbb{CP}^1$ (by forgetting one of the factors), which exhibit $X$ as a ramified double cover of $\mathbb{CP}^1\times\mathbb{CP}^1$. The three covering involutions are now denoted by $\sigma_1,\sigma_2,\sigma_3$, and $T=\sigma_1\circ\sigma_2\circ\sigma_3$ is an automorphism of $X$ with $h(T)=\log(9+4\sqrt{5})\sim 2.88$ (see \cite{Can}). Together they generate a subgroup of $\mathrm{Aut}(X)$ isomorphic to the free product $\mathbb{Z}/2\mathbb{Z} * \mathbb{Z}/2\mathbb{Z} * \mathbb{Z}/2\mathbb{Z}.$
\end{example}

\begin{example}[McMullen \cite{McM}]\label{mc} Using the refined Torelli Theorem, and substantial work, McMullen has constructed \cite{McM} examples of non-projective $K3$ surfaces $X$ with automorphisms $T$ with $h(T)>0$ which admit a Siegel disc: this is an open subset $\Delta\subset X$ preserved by $T$ and biholomorphic to the bidisc in $\mathbb{C}^2$, such that in $\Delta$ the automorphism $T$ is holomorphically conjugate to an irrational rotation $(z_1,z_2)\mapsto (a z_1,bz_2)$ of the bidisc (which means that $|a|=|b|=1$ and $a,b$ and $ab$ are not roots of unity).
\end{example}

It is also worth remarking here that a result of Cantat \cite{Can2} shows that the only compact complex surfaces which admit automorphisms with positive topological entropy are $K3$, Enriques, $2$-tori, iterated blowups of these, and blowups of $\mathbb{CP}^2$ at $k$ points with $k\geq 10$. The dynamical study of such automorphisms on $2$-tori is elementary, on Enriques surfaces it can be reduced to the $K3$ surface that is its double cover, and on blowups of $K3$, Enriques and tori it can be reduced to the base case. Thus, the only ``interesting'' cases are $K3$ surfaces and blowups of $\mathbb{CP}^2$. See e.g. \cite{CaT,McM3} and references therein for more on automorphisms of rational surfaces.

\subsection{Eigenclasses}
From now on we assume that $T:X\to X$ is a $K3$ automorphism with $h(T)>0$. Since $T^*:\mathcal{H}\to\mathcal{H}$ is a hyperbolic isometry, it has two fixed points on its ideal boundary $\de\mathcal{H}$, which correspond to two nontrivial eigenclasses $[\eta_+]$ and $[\eta_-]\in H^{1,1}(X,\mathbb{R})$ which satisfy
$$T^*[\eta_{\pm}]=e^{\pm h}[\eta_{\pm}], \quad\int_X[\eta_{\pm}]^2=0,$$
and up to rescaling these classes if necessary we may assume also that
$$\int_X[\eta_+]\wedge[\eta_-]=1.$$
Using the Gromov-Yomdin Theorem $h=h(T)=\log\rho$, it follows easily that for any given K\"ahler metric $\omega$ on $X$,
$$\lim_{n\to+\infty}\frac{(T^n)^*[\omega]}{e^{nh}}=\left(\int_X[\omega]\wedge[\eta_-]\right)[\eta_+], \quad \lim_{n\to+\infty}\frac{(T^{-n})^*[\omega]}{e^{nh}}=\left(\int_X[\omega]\wedge[\eta_+]\right)[\eta_-].$$
This implies that the classes $[\eta_\pm]$ belong to $\de\mathcal{C}_X$, where recall that $\mathcal{C}_X$ is the cone in $H^{1,1}(X,\mathbb{R})$ of cohomology classes of the form $[\omega]$ for some K\"ahler metric $\omega$ on $X$. These classes are irrational in a strong sense, namely that the line $\mathbb{R}.[\eta_+]$ intersects $H^2(X,\mathbb{Q})$ only in the origin (and the same holds for $[\eta_-]$).

\subsection{Eigencurrents}
Let us now look in more detail at the eigenclasses $[\eta_{\pm}]$. Fix two closed real $(1,1)$-forms $\alpha_+$ and $\alpha_-$ with $[\alpha_\pm]=[\eta_\pm]$.
Every closed real $(1,1)$-form $\beta$ in the cohomology class $[\eta_+]$ can therefore be written as $\beta=\alpha_++\ddbar\vp$ for some smooth function $\vp$. We will write $\beta\geq 0$ if $\beta$ is Hermitian semipositive at every point. Given that the class $[\eta_+]$ is on the boundary of the K\"ahler cone $\mathcal{C}_X$, and that classes inside this cone contain representatives which are smooth and strictly positive, one might naively expect that the classes $[\eta_\pm]$ might contain smooth semipositive representatives. This is in general not the case, as we shall see below (cf. Theorem \ref{eig}), however it is possible to find a semipositive representative if one is willing to relax smoothness.

More precisely, a closed positive $(1,1)$-current $\beta$ in the class $[\eta_+]$ is a $(1,1)$-form with distributional coefficients which can be written as $\beta=\alpha_++\ddbar\vp$ where $\vp$ is quasi-psh (i.e. in local charts it equals the sum of a plurisubharmonic function and a smooth function), such that $T\geq 0$ holds in the weak sense (which means that $\langle \beta,i\xi\wedge\ov{\xi}\rangle\geq 0$ for all smooth $(1,0)$-forms $\xi$). Such currents can be pulled back via $T$ by defining $T^*\beta=T^*\alpha_++\ddbar(\vp\circ T)$, and the pullback is a closed positive $(1,1)$-current in the class $[T^*\eta_+]$.
We then have the following crucial result:

\begin{theorem}[Cantat \cite{Can}]\label{eig}
Let $T:X\to X$ be an automorphism of a $K3$ surface with $h(T)>0$. Then
\begin{itemize}
\item[(a)] The classes $[\eta_\pm]$ contain a unique closed positive $(1,1)$-current $\eta_\pm$
\item[(b)] These currents satisfy $T^*\eta_\pm=e^{\pm h}\eta_{\pm}$
\item[(c)] (Dinh-Sibony \cite{DS}) We can write $\eta_\pm=\alpha_\pm+\ddbar\vp_\pm$ where $\vp_\pm$ is quasi-psh and $C^\alpha(X)$ for some $\alpha>0$
\item[(d)] The wedge product $\mu=\eta_+\wedge\eta_-$ exists by Bedford-Taylor \cite{BT} and (c), and is a $T$-invariant probability measure on $X$
\item[(e)] $\mu$ is mixing, hence ergodic, and is the unique measure of maximal entropy
\end{itemize}
\end{theorem}
Let us give a few explanations about this result. About part (c), the fact that $\vp_\pm$ is continuous was proved earlier in \cite{Can} when $X$ is projective, see also \cite{DG} for a concise exposition of the H\"older regularity result of Dinh-Sibony. In part (d), the wedge product of Bedford-Taylor \cite{BT} is defined as $\mu=\eta_+\wedge\alpha_-+\ddbar(\vp_-\eta_+)$, which is well-defined since the distributional coefficients of $\eta_+$ are in fact measures and $\vp_-$ is H\"older. About part (e), $\mu$ being mixing means that for every $f,g$ $\mu$-measurable functions,
$$\int_X f(T^n(x)) g(x)d\mu(x)\overset{n\to+\infty}{\to}\int_X fd\mu\int_X gd\mu,$$
which implies that $\mu$ is ergodic (every $T$-invariant $\mu$-measurable subset of $X$ has either zero or full measure). Lastly, for every $T$-invariant (Borel) probability measure $\nu$ on $X$, one has (see e.g. \cite{Ne}) the Kolmogorov-Sinai entropy $h_\nu(T)$ of $\nu$, which is always bounded above by the topological entropy $h(T)$. If $h_\nu(T)=h(T)$, then $\nu$ is called a measure of maximal entropy. Ergodic measures of maximal entropy always exist in our setting by a general result of Newhouse \cite{Ne}, and part (e) then asserts that there is only one such measure, $\mu$. Part (e) was proved in \cite{Can} for projective $K3$ surfaces, and in \cite{DD} in general.

\section{From geometry to dynamics: Kummer rigidity}\label{s1}
\subsection{Two invariant measures}As discussed in the previous section, for an automorphism $T:X\to X$ of a $K3$ surface with positive topological entropy, one obtains two natural $T$-invariant probability measures on $X$, the ``Lebesgue'' measure $\dVol=\Omega\wedge\ov{\Omega}$ (it is in fact equal to the Lebesgue measures in suitable local coordinate charts) and the measure $\mu$ of maximal entropy. It is then natural to ask about the relation between them.

We see from Example \ref{mc} that in general $\mu\neq \dVol$, since $\mu$ is ergodic but in Example \ref{mc} the Lebesgue measure cannot be ergodic since $T$ is a rotation on the Siegel disc. On the other hand, it is not hard to check that if $(X,T)$ is a Kummer example (see Example \ref{ku}) then we do indeed have $\mu=\dVol$.

\subsection{Kummer rigidity} Cantat \cite[p.162]{CaT} and McMullen \cite[Conjecture 3.31]{McM2} had conjectured the following:

\begin{conjecture}\label{rig}
Let $T:X\to X$ be a $K3$ automorphism with positive topological entropy. Then $\mu\ll \dVol$ if and only if $(X,T)$ is a Kummer example.
\end{conjecture}

In other words, $\mu$ being absolutely continuous with respect to Lebesgue suffices to conclude that $(X,T)$ is a Kummer example, and then {\em a posteriori} $\mu=\dVol$. This ``Kummer rigidity'' conjecture is analogous to a rigidity theorem for rational maps of $\mathbb{CP}^1$ of Zdunik \cite{Zd}, and for general endomorphisms of $\mathbb{CP}^n$ in \cite{BD,BL}, where the role of Kummer example is played by Latt\`es maps.

Furthermore, McMullen \cite{McM2} also formulated an extension of this conjecture to automorphisms with positive topological entropy of general compact complex surfaces (which are necessarily K\"ahler), where the unique measure of maximal entropy $\mu$ still exists \cite{DD}, and he again conjectured that $\mu$ is absolutely continuous with respect to the Lebesgue measure if and only if $(X,T)$ is a generalized Kummer example (which need not be a $K3$ surface). Since this survey focuses on $K3$ surfaces, we refer the reader to \cite{CD} for more details.

Conjecture \ref{rig}, together with the more general one for other surfaces, was settled by Cantat-Dupont \cite{CD} when $X$ is projective. In \cite{FT2}, Filip and the author gave a different proof for $K3$ surfaces, which does not require projectivity, and uses the Ricci-flat metrics:

\begin{theorem}[Filip-T. \cite{FT2}, Cantat-Dupont \cite{CD} when $X$ projective]\label{rigg}Conjecture \ref{rig} is true, and in fact the following are equivalent for $K3$ automorphisms with positive topological entropy:
\begin{itemize}
\item[(a)]$\mu\ll\dVol$
\item[(b)]$\mu=\dVol$
\item[(c)]The eigencurrents $\eta_\pm$ are smooth (or just continuous off a closed analytic subset)
\item[(d)]$(X,T)$ is a Kummer example
\end{itemize}
\end{theorem}

Combining the results in \cite{CD} and \cite{FT2}, one also obtains the proof of the more general conjecture for arbitrary surfaces, since projective surfaces are covered in \cite{CD}, and the only non-projective ones which need to be dealt with are $K3$ for which \cite{FT2} applies.

As a corollary of Theorem \ref{rigg}, it follows that in Example \ref{mc}, the measure $\mu$ cannot be absolutely continuous with respect to Lebesgue, since as we remarked earlier $\mu\neq\dVol$. In this case in fact it is easy to see that $\eta_\pm$ (and so also $\mu$) vanish on the Siegel disc (see \cite[Theorem 11.4]{McM}).

\subsection{Ricci-flat metrics and rigidity}\label{sr}Let us give a sketch of proof of (part of) Theorem \ref{rigg}. It is easy to see that $(d)\Rightarrow (c)\Rightarrow(b)\Rightarrow(a)$, and we will discuss the proof that $(b)\Rightarrow(d)$ under the extra assumption that $X$ contains no $T$-periodic curves. When $X$ contains periodic curves, these can be contracted to obtain an orbifold $K3$ surface with an induced automorphism with the same positive entropy, and the arguments we will describe have to be applied on the orbifold. We will not discuss here the implication $(a)\Rightarrow(b)$, which involves other ingredients.

To start, let us look at the cohomology class $[\eta_+]+[\eta_-]$. It belongs to the closure of the K\"ahler cone, and it satisfies
\begin{equation}\label{bs}
\int_X([\eta_+]+[\eta_-])^2=2\int_X[\eta_+]\wedge[\eta_-]=2,
\end{equation}
so it is a {\em nef and big} class using terminology borrowed from algebraic geometry. Consider the {\em null locus} of $[\eta_+]+[\eta_-]$,
$$\mathrm{Null}([\eta_+]+[\eta_-])=\bigcup_{\int_C([\eta_+]+[\eta_-])=0}C,$$
where the union is over all irreducible (complex) curves $C\subset X$ with $\int_C([\eta_+]+[\eta_-])=0$. By general results of Collins and the author \cite{CT} (which hold for nef and big classes on arbitrary compact K\"ahler manifolds) the null locus is a closed analytic subset of $X$, and so it consists of the union of finitely many irreducible curves.

In fact, $\mathrm{Null}([\eta_+]+[\eta_-])$ is the same as the union of the $T$-periodic curves. Indeed, if $C\subset X$ is any irreducible curve, we have
\begin{equation}\label{eq}
\int_C[\eta_\pm]=\int_{T^{-N}(C)}(T^N)^*[\eta_\pm]=e^{\pm Nh}\int_{T^{-N}(C)}[\eta_\pm],
\end{equation}
for all $N\in\mathbb{Z}$, from which it follows that if $C$ is periodic ($T^{-N}(C)=C$ for some $N\in\mathbb{Z}\backslash\{0\}$) then $\int_C[\eta_\pm]=0$ so $C\subset \mathrm{Null}([\eta_+]+[\eta_-])$. Conversely, if $C\subset \mathrm{Null}([\eta_+]+[\eta_-])$, then $\int_{C}([\eta_+]+[\eta_-])=0$ and $\int_C[\eta_\pm]\geq 0$ since $[\eta_\pm]$ belong to the closure of $\mathcal{C}_X$, and so $\int_C[\eta_\pm]=0$, and from \eqref{eq} we get
$\int_{T^{-N}(C)}[\eta_\pm]=0$ for all $N\in\mathbb{Z}$, so all the irreducible curves $T^{-N}(C)$ are contained in $\mathrm{Null}([\eta_+]+[\eta_-])$. Since this consists of finitely many irreducible curves, we get $T^{-N}(C)=T^{-M}(C)$ for some distinct integers $N,M$, and so $C$ is $T$-periodic.

Going back to our main argument, since we assume that there are no $T$-periodic curves, we have $\mathrm{Null}([\eta_+]+[\eta_-])=\emptyset$, which e.g. by \cite{CT} implies that the class $[\eta_+]+[\eta_-]$ is K\"ahler, and by Yau's Theorem \ref{yau} we can fix a Ricci-flat K\"ahler metric $\omega$ on $X$ in this class. Thanks to \eqref{bs}, it satisfies
$$\omega^2 =2\dVol.$$
Then for $N\geq 1$, we let $$\omega_N=(T^N)^*\omega,$$
which is the unique Ricci-flat K\"ahler metric in the class $$(T^N)^*([\eta_+]+[\eta_-])=e^{Nh}[\eta_+]+e^{-Nh}[\eta_-],$$ and also satisfies
$$\omega_N^2=2\dVol.$$
For each $N\geq 1$, define now a function $\lambda(x,N)$, which is continuous in $x$, so that the largest eigenvalue of $\omega_N(x)$ with respect to $\omega(x)$ is equal to $e^{2\lambda(x,N)}$. Since $\omega^2=\omega_N^2$, it follows that $\lambda(x,N)\geq 0$ and that the smallest eigenvalue of  $\omega_N(x)$ with respect to $\omega(x)$ is equal to $e^{-2\lambda(x,N)}$, and so the trace of $\omega_N(x)$ with respect to $\omega(x)$ equals
$$\tr{\omega}{\omega_N}(x)=2\frac{\omega\wedge\omega_N}{\omega^2}(x)=e^{2\lambda(x,N)}+e^{-2\lambda(x,N)}.$$
Before we continue with our arguments, we need the following crucial claim:
\begin{equation}\label{ly}
2\int_X \lambda(x,N)\dVol(x)\geq Nh,
\end{equation}
for all $N\geq 1$. To see this, first note that $\lambda(x,N)=\log\| D_xT^N\|_\omega$ from which a standard argument (see \cite[\S 2.2]{FT2}) shows that $I_N=\int_X \lambda(x,N)\dVol(x)$ is subadditive and so $\Lambda=\lim_{N\to+\infty}\frac{I_N}{N}$ exists and satisfies
\begin{equation}\label{ly2}
\Lambda\leq \frac{I_N}{N},
\end{equation}
for all $N$. The number $\Lambda$ is in fact the largest {\em Lyapunov exponent} of $\dVol$, since we assume that $\dVol=\mu$, and we know that in general $\mu$ is ergodic (Theorem \ref{eig} (e)). The Ledrappier-Young formula \cite{LY2} then gives that the topological entropy $h$, which equals the Kolmogorov-Sinai entropy of $\mu=\dVol$ (recall again Theorem \ref{eig} (e)), is equal to
\begin{equation}\label{ly3}
h=\Lambda\cdot \dim_+(\mu)=2\Lambda,
\end{equation}
since $\dim_+(\mu)$ (the dimension of $\mu$ along the unstable directions) equals $2$ because $\mu=\dVol$. Combining \eqref{ly2} and \eqref{ly3} gives \eqref{ly}.

We now use \eqref{ly} for our main computation as follows: by Stokes's Theorem, the integral $\int_X\omega\wedge\omega_N$ can be compute in cohomology as
$$\int_X\omega\wedge\omega_N=\int_X([\eta_+]+[\eta_-])\wedge(e^{Nh}[\eta_+]+e^{-Nh}[\eta_-])=e^{Nh}+e^{-Nh},$$
and so using Jensen's inequality
\begin{equation}\label{1}\begin{split}
\log(e^{Nh}+e^{-Nh})&=\log\left(\int_X\omega\wedge\omega_N\right)\geq\int_X\log\left(\frac{\omega\wedge\omega_N}{\dVol}\right)\dVol\\
&=\int_X\log\left(\frac{2\omega\wedge\omega_N}{\omega^2}\right)\dVol\\
&=\int_X\log\left(e^{2\lambda(x,N)}+e^{-2\lambda(x,N)}\right)\dVol(x),
\end{split}\end{equation}
but noting that $t\mapsto\log(e^t+e^{-t})$ is convex and increasing for $t\geq 0$, we can apply Jensen's inequality again and \eqref{ly} to get
\begin{equation}\label{2}\begin{split}
\int_X\log\left(e^{2\lambda(x,N)}+e^{-2\lambda(x,N)}\right)\dVol(x)&\geq
\log\left(e^{2\int_X\lambda(x,N)\dVol(x)}+e^{-2\int_X\lambda(x,N)\dVol(x)}\right)\\
&\geq \log(e^{Nh}+e^{-Nh}),
\end{split}\end{equation}
which implies that all the inequalities in \eqref{1} and \eqref{2} must be equalities and so $\lambda(x,N)=\frac{Nh}{2}$ holds for all $x\in X$ and $N\geq 1$. Going back to the definition of $\lambda(x,N)$, this means that we obtain two $\omega$-orthogonal $T$-invariant line subbundles of $TX$, one expanded and one contracted by $T$. By Ghys \cite[Proposition 2.2]{Gh}, these give two holomorphic foliations on $X$ which are preserved by $T$, which is already enough to conclude that $(X,T)$ is a Kummer example by Cantat \cite[Theorem 7.4]{Can} (or \cite[Theorem 3.1]{CF} which only needs one invariant foliation). Alternatively, one can directly use these two invariant foliations to show that $\omega$ must be flat, and then that $(X,T)$ is Kummer, see \cite[\S 3.2]{FT2}. This concludes our sketch of the proof that  $(b)\Rightarrow(d)$ in Theorem \ref{rigg}.
\section{From dynamics to geometry: limits of Ricci-flat metrics}\label{s2}
In the previous section we saw an application of the Ricci-flat K\"ahler metrics on $K3$ surfaces to dynamics. Here we go in the opposite direction, and use dynamics to prove results about the Ricci-flat metrics.

Let $X$ be a $K3$ surface. Recall that thanks to Yau's Theorem \ref{yau} for every K\"ahler class $[\alpha]\in \mathcal{C}_X\subset H^{1,1}(X,\mathbb{R})$ (an open convex cone in this cohomology group) there is a unique Ricci-flat K\"ahler metric $\omega$ with $[\omega]=[\alpha]$. A natural question to ask is how do these metrics behave as the class $[\alpha]$ varies. It is easy to see (either from Yau's explicit estimates, or using the implicit function theorem) that the Ricci-flat metrics vary continuously in the smooth topology as long as their cohomology class is contained in a fixed relatively compact subset of $\mathcal{C}_X$. We would then like to know what happens when we approach a limiting class $[\alpha]\in \de\mathcal{C}_X$.

This is a problem that has received much attention recently, see for example the author's surveys \cite{To2,To3,To4} and references therein. We will just focus on the following basic setup: given a class $[\alpha]\in\de\mathcal{C}_X$, and a fixed Ricci-flat K\"ahler metric $\omega$ on $X$, let $\omega_t, 0<t\leq 1$ be the unique Ricci-flat K\"ahler metric on $X$ cohomologous to $[\alpha]+t[\omega]$. What is the behavior of $\omega_t$ as $t\to 0$?

\subsection{Noncollapsed limits}Suppose first that $\int_X[\alpha]^2>0$. In this case, as discussed earlier, the null locus $\mathrm{Null}([\alpha])$ of $[\alpha]$, which is the union of all irreducible curves which intersect trivially with $[\alpha]$, is a closed analytic subset of $X$. Then, as shown in \cite{To0} and \cite{CT}, the Ricci-flat metrics $\omega_t$ converge locally smoothly (as tensors) on compact sets away from $\mathrm{Null}([\alpha])$ to a Ricci-flat K\"ahler metric $\omega_0$ on $X\backslash \mathrm{Null}([\alpha])$. See \cite{To4} for more information and higher-dimensional generalizations.

\subsection{Collapsed fibrations limits}Suppose next that $\int_X[\alpha]^2=0$ and that $[\alpha]=\pi^*[\omega_{\mathbb{CP}^1}]$ is the pullback of a K\"ahler class on $\mathbb{CP}^1$ via an elliptic fibration $\pi:X\to\mathbb{CP}^1$. Then, as shown in \cite{GW} when $\pi$ has $24$ singular fibers of type $I_1$ and in \cite{GTZ,HT} in general, the Ricci-flat metrics $\omega_t$ again converge locally smoothly (as tensors) on compact sets away from the singular fibers $S\subset X$ (a closed analytic subset of $X$) to the pullback of a K\"ahler metric $\omega_0$ on $\mathbb{CP}^1\backslash \pi(S)$. The limiting metric $\omega_0$ is not Ricci-flat, its Ricci curvature is a Weil-Petersson semipositive form that measures the variation of complex structure of the smooth fibers. See again \cite{To4} for more details and generalizations.

It is also interesting to note that if $0\neq [\alpha]\in\de\mathcal{C}_X$ satisfies $\int_X[\alpha]^2=0$ and $[\alpha]\in H^2(X,\mathbb{Q})$, then in fact $[\alpha]=\pi^*[\omega_{\mathbb{CP}^1}]$ for some elliptic fibration on $X$ (see \cite[Proposition 1.4]{FT}).
\subsection{Enter dynamics}Based on the two previous results, the author had conjectured in \cite{To2,To3} that for arbitrary classes $[\alpha]\in\de\mathcal{C}_X$, the Ricci-flat metrics $\omega_t$ should converge locally smoothly on compact sets away from some closed analytic subset $S$ of $X$. However, this turns out to be false, as observed by Filip and the author \cite{FT}:

\begin{theorem}\label{st}
Let $X$ be a $K3$ surface with an automorphism $T$ with positive topological entropy such that $(X,T)$ is not a Kummer example (for example, those described in Examples \ref{we}, \ref{222} and \ref{mc}), let $[\alpha]=[\eta_+]$ and $\omega_t$ the Ricci-flat K\"ahler metric on $X$ cohomologous to $[\eta_+]+t[\omega]$, $0<t\leq 1$. Then as $t\to 0$ the metrics $\omega_t$ cannot converge in $C^0_{\rm loc}$ on the complement of any closed analytic subset of $X$.
\end{theorem}

Indeed, this is essentially a corollary of Theorem \ref{rigg}: by weak compactness of currents, it is easy to show that the metrics $\omega_t$ must converge in the weak topology of currents to the eigencurrent $\eta_+$ as $t\to 0$ (here we use that $\eta_+$ is the unique closed positive current in its class by Theorem \ref{eig} (a)), so if $\omega_t$ was also converging in $C^0_{\rm loc}(X\backslash S)$ for some closed analytic subset $S$, then $\eta_+$ would be continuous on $X\backslash S$. Now, if this was true for both $\eta_+$ and $\eta_-$, then Theorem \ref{rigg} would immediately give a contradiction (so Theorem \ref{st} follows if we allow perhaps replacing $[\eta_+]$ by $[\eta_-]$). To show that just continuity of $\eta_+$ on $X\backslash S$ is enough to conclude that $\mu\ll\dVol$ (and hence derive a contradiction by Theorem \ref{rigg} again) one needs to work just a little bit more, using \cite{DD,LY}, see \cite[Theorem 3.3 (3)]{FT}.

\subsection{Other boundary classes}
To conclude, we discuss briefly what is expected to happen to the Ricci-flat K\"ahler metrics $\omega_t$ when their cohomology class approaches $0\neq [\alpha]\in\de\mathcal{C}_X$ which satisfies $\int_X[\alpha]^2=0$ but does not come from the base of an elliptic fibration, and is not an eigenclass for an automorphism with positive entropy.

We fix a smooth representative $\alpha$ of its class, which is a closed real $(1,1)$-form. Since $[\alpha]$ is a limit of K\"ahler classes, weak compactness of currents easily shows that there are closed positive $(1,1)$-currents $\beta=\alpha+\ddbar\vp_0$ in the class $[\alpha]$, with $\vp_0$ quasi-psh normalized by $\sup_X\vp_0=0$ say. This is again expected to be unique, see the ongoing work of Sibony-Verbitsky \cite{SV}.

The Ricci-flat metrics $\omega_t$ in the class $[\alpha]+t[\omega], 0<t\leq 1$, can be written as $\omega_t=\alpha+t\omega+\ddbar\vp_t$, where $\vp_t$ are smooth functions which are uniquely determined if we normalize them by $\sup_X\vp_t=0$.

\begin{conjecture}
Let $X$ be a $K3$ surface, $\alpha$ a closed real $(1,1)$-form with $0\neq[\alpha]\in\de\mathcal{C}_X$ and $\int_X\alpha^2=0$. Let $\omega$ be a K\"ahler metric on $X$, and for $0<t\leq 1$ let $\omega_t=\alpha+t\omega+\ddbar\vp_t$ be the Ricci-flat K\"ahler metric in the class $[\alpha]+t[\omega]$ with normalization $\sup_X\vp_t=0$. Then there is a closed positive $(1,1)$-current $\beta=\alpha+\ddbar\vp_0\geq 0$ with $\vp_0\in C^0(X), \sup_X\vp_0=0$, and
\begin{equation}\label{conv}
\vp_t\to\vp_0,
\end{equation}
uniformly on $X$ as $t\to 0$.
\end{conjecture}

This conjecture is known in the case when $[\alpha]$ is an eigenclass for an automorphism with positive entropy, since in this case $\vp_0$ is even $\gamma$-H\"older continuous for some $\gamma>0$ by Theorem \ref{eig} (c), and the convergence of $\vp_t$ to $\vp_0$ is easily seen to hold in $C^\gamma(X)$.

Interestingly, this conjecture is not known when $[\alpha]$ comes from the base of an elliptic fibration $\pi:X\to\mathbb{CP}^1$: in this case we do know that $\vp_0\in C^\gamma(\mathbb{CP}^1)$ for some $\gamma>0$ (since its Laplacian is globally in $L^p$ for some $p>1$ \cite[Corollary 3.1]{ST0}), but the global convergence in \eqref{conv} uniformly on all of $X$ (not just away from the singular fibers) is unknown.

And of course the most interesting case is when $[\alpha]$ is neither an eigenclass nor comes from an elliptic fibration, in which case even the existence of a continuous $\vp_0$ as above is unknown.

\section{Some conjectures}\label{sc}
In this last section we briefly discuss a few open problems related to the dynamics of automorphisms of $K3$ surfaces that the author learned from S. Filip, see also Cantat's ICM paper \cite{CanS3} for many other problems.

\subsection{Positive Lyapunov exponent}
Let $T:X\to X$ be a $K3$ automorphism with positive topological entropy $h>0$ and fix a Ricci-flat K\"ahler metric $\omega$ on $X$. The largest Lyapunov exponent of $\dVol$ (which appeared in section \ref{sr} in the special case when $\dVol=\mu$) is then defined as
$$\Lambda=\int_X\left(\lim_{N\to+\infty}\frac{h}{2N}\log\| D_xT^N\|_\omega\right)\dVol(x).$$
This is easily seen to be finite, and if we let $\omega_N=(T^N)^*\omega$ then $\log\| D_xT^N\|_\omega$ is equal to the quantity $\lambda(x,N)\geq 0$ defined in section \ref{sr} (namely the largest eigenvalue of $\omega_N(x)$ with respect to $\omega(x)$ is $e^{2\lambda(x,N)}$). This shows that $\Lambda \geq 0$, and the major outstanding problem is then (see also the discussion in Cantat's thesis \cite[Chapter 3]{CaT}):
\begin{conjecture}\label{c2} Let $T:X\to X$ be an automorphism with positive topological entropy of a projective $K3$ surface. Then we have $\Lambda>0$. 
\end{conjecture}
It would already be extremely interesting to show that in the setting of Conjecture \ref{c2} there is dense $T$-orbit. Furthermore, once $\Lambda>0$ one expects more:
\begin{question}\label{c1} Let $T:X\to X$ be a $K3$ automorphism with positive topological entropy and suppose that $\Lambda>0$. Does it follow that $\dVol$ is $T$-ergodic?
\end{question}
Recall that the measure of maximal entropy $\mu$ is always $T$-ergodic (even mixing), and it also has positive Lyapunov exponent by the Ledrappier-Young formula \cite{LY2}, but in general $\mu$ is quite different from $\dVol$ as shown in Theorem \ref{rigg}. 
\subsection{The support of $\mu$}
Let again $T:X\to X$ be a $K3$ automorphism with positive topological entropy $h>0$, and let $\mu=\eta_+\wedge\eta_-$ be the measure with maximal entropy from Theorem \ref{eig}. By Theorem \ref{rigg} we know that if $(X,T)$ is not a Kummer example then $\mathrm{Supp}(\mu)$ is a Lebesgue null-set. Nevertheless, this set should be quite fractal, and we expect that:
\begin{conjecture}\label{c3}Let $T:X\to X$ be an automorphism with positive topological entropy of a projective $K3$ surface. Then $\ov{\mathrm{Supp}(\mu)}$ has full Lebesgue measure.
\end{conjecture}
Thanks to a result of Dinh-Sibony (see \cite[Theorem 7.6]{CanS2}), an affirmative answer to this conjecture would give a negative answer to \cite[Question 3.4]{CanS3}. 
Note that this conjecture is false when $X$ is not projective, as shown by McMullen's examples of $K3$ automorphisms with Siegel discs in Example \ref{mc} (which are not projective): indeed, $\mu$ vanishes completely on the Siegel disc. It seems quite likely that in fact the Siegel disc, when it exists, is rather large:
\begin{question}\label{c4}Let $T:X\to X$ be an automorphism with positive topological entropy of a $K3$ surface which admits a Siegel disc $\Delta\subset X$. What is the Lebesgue measure of $\ov{\mathrm{Supp}(\mu)}$? Could it be zero? 
\end{question}
In other words, are there Siegel discs so that the Lebesgue measure of $X\backslash \Delta$ is zero?


\begin{thebibliography}{99}
\bibitem{BHPV}  W.P. Barth, K. Hulek, C.A.M. Peters, A. Van de Ven, {\em Compact complex surfaces. Second edition}, Ergebnisse der Mathematik und ihrer Grenzgebiete. 3. Folge. 4. Springer-Verlag, Berlin, 2004.
\bibitem{Be} A. Beauville, {\em Vari\'et\'es K\"ahleriennes dont la premi\`ere classe de Chern est nulle}, J. Differential Geom. {\bf 18} (1983), no. 4, 755--782.
\bibitem{BT} E. Bedford, B.A. Taylor, {\em The Dirichlet problem for a complex Monge-Amp\`ere equation}, Invent. Math. {\bf 37} (1976), no. 1, 1--44.
\bibitem{BD} F. Berteloot, C. Dupont, {\em Une caract\'erisation des endomorphismes de Latt\`es par leur mesure de Green}, Comment. Math. Helv. {\bf 80} (2005), no. 2, 433--454.
\bibitem{BL} F. Berteloot, J.-J. Loeb, {\em Une caract\'erisation g\'eom\'etrique des exemples de Latt\`es de $\mathbb{P}^k$}, Bull. Soc. Math. France {\bf 129} (2001), no. 2, 175--188.
\bibitem{Bu} N. Buchdahl, {\em On compact K\"ahler surfaces}, Ann. Inst. Fourier (Grenoble) {\bf 49} (1999), no. 1, 287--302.
\bibitem{BR} D. Burns, M. Rapoport, {\em On the Torelli problem for k\"ahlerian $K--3$ surfaces}, Ann. Sci. \'Ecole Norm. Sup. (4) {\bf 8} (1975), no. 2, 235--273.
\bibitem{Cal} E. Calabi, {\em On K\"ahler manifolds with vanishing canonical class}, in {\em Algebraic geometry and topology. A symposium in honor of S. Lefschetz},  pp. 78--89. Princeton University Press, Princeton, N. J., 1957.
\bibitem{CaT} S. Cantat, {\em Dynamique des automorphismes des surfaces complexes compactes}, PhD Thesis, \'Ecole Normale Sup\'erieure de Lyon, 1999.
\bibitem{Can2} S. Cantat, {\em Dynamique des automorphismes des surfaces projectives complexes}, C. R. Acad. Sci. Paris S\'er. I Math. {\bf 328} (1999), no. 10, 901--906.
\bibitem{Can} S. Cantat, {\em Dynamique des automorphismes des surfaces $K3$}, Acta Math. {\bf 187} (2001), no. 1, 1--57.
\bibitem{CanS} S. Cantat, {\em Quelques aspects des syst\`emes dynamiques polynomiaux: existence, exemples, rigidit\'e}, in {\em Quelques aspects des syst\`emes dynamiques polynomiaux}, 13--95, Panor. Synth\`eses, 30, Soc. Math. France, Paris, 2010.
\bibitem{CanS2} S. Cantat, {\em Dynamics of automorphisms of compact complex surfaces}, in {\em Frontiers in complex dynamics}, 463--514, Princeton Math. Ser., 51, Princeton Univ. Press, Princeton, NJ, 2014.
\bibitem{CanS3} S. Cantat, {\em Automorphisms and dynamics: a list of open problems}, in {\em Proceedings of the International Congress of Mathematicians--Rio de Janeiro 2018. Vol. II. Invited lectures}, 619--634, World Sci. Publ., Hackensack, NJ, 2018.
\bibitem{CD} S. Cantat, C. Dupont, {\em Automorphisms of surfaces: Kummer rigidity and measure of maximal entropy}, J. Eur. Math. Soc. (JEMS) {\bf 22} (2020), no. 4, 1289--1351.
\bibitem{CF} S. Cantat, C. Favre, {\em Sym\'etries birationnelles des surfaces feuillet\'ees}, J. Reine Angew. Math. {\bf 561} (2003), 199--235; Corrigendum, ibid. {\bf 582} (2005), 229--231.
\bibitem{CT} T. Collins, V. Tosatti, {\em K\"ahler currents and null loci}, Invent. Math. {\bf 202} (2015), no.3, 1167--1198.
\bibitem{DD} H. De Th\'elin, T.-C. Dinh, {\em Dynamics of automorphisms on compact K\"ahler manifolds}, Adv. Math. {\bf 229} (2012), no. 5, 2640--2655.
\bibitem{DG} J. Diller, V. Guedj, {\em Regularity of dynamical Green's functions}, Trans. Amer. Math. Soc. {\bf 361} (2009), no. 9, 4783--4805.
\bibitem{DS} T.-C. Dinh, N. Sibony, {\em Green currents for holomorphic automorphisms of compact K\"ahler manifolds}, J. Amer. Math. Soc. {\bf 18} (2005), no. 2, 291--312.
\bibitem{Fil} S. Filip, {\em An introduction to $K3$ surfaces and their dynamics}, lecture notes. \url{http://math.uchicago.edu/~sfilip/public_files/lectures_k3_dynamics.pdf}
\bibitem{FT} S. Filip, V. Tosatti, {\em Smooth and rough positive currents}, Ann. Inst. Fourier (Grenoble) {\bf 68} (2018), no.7, 2981--2999.
\bibitem{FT2} S. Filip, V. Tosatti, {\em Kummer rigidity for $K3$ surface automorphisms via Ricci-flat metrics}, to appear in Amer. J. Math.
\bibitem{K3} {\em G\'eom\'etrie des surfaces $K3$: modules et p\'eriodes}, Papers from the seminar held in Palaiseau, October 1981-January 1982. Ast\'erisque No. {\bf 126} (1985). Soci\'et\'e Math\'ematique de France, Paris, 1985. pp. 1--193.
\bibitem{Gh} \'E. Ghys, {\em Holomorphic Anosov systems}, Invent. Math. {\bf 119} (1995), no. 3, 585--614.
\bibitem{Gr} M. Gromov, {\em On the entropy of holomorphic maps}, Enseign. Math. (2) {\bf 49} (2003), no. 3-4, 217--235.
\bibitem{GTZ} M. Gross, V. Tosatti, Y. Zhang, {\em Collapsing of abelian fibred Calabi-Yau manifolds}, Duke Math. J. {\bf 162} (2013), no. 3, 517--551.
\bibitem{GW} M. Gross, P.M.H. Wilson, \emph{Large complex structure limits of $K3$ surfaces}, J. Differ. Geom. {\bf 55} (2000), no. 3, 475--546.
\bibitem{HT} H.-J. Hein, V. Tosatti, {\em Remarks on the collapsing of torus fibered Calabi-Yau manifolds},  Bull. Lond. Math. Soc. {\bf 47} (2015), no. 6, 1021--1027.
\bibitem{Huy} D. Huybrechts, {\em Lectures on $K3$ surfaces}, Cambridge Studies in Advanced Mathematics, 158. Cambridge University Press, Cambridge, 2016.
\bibitem{Ko} K. Kodaira, {\em On the structure of compact complex analytic surfaces, I}, Amer. J. Math. {\bf 86} (1964), no. 4, 751--798.
\bibitem{Kon} S. Kond$\overline{\mbox{o}}$, {\em $K3$ surfaces}, EMS Tracts in Mathematics, 32. European Mathematical Society (EMS), Z\"urich, 2020.
\bibitem{Ku} V.S. Kulikov, {\em Degenerations of $K3$ surfaces and Enriques surfaces}, Izv. Akad. Nauk SSSR Ser. Mat. {\bf 41} (1977), no. 5, 1008--1042.
\bibitem{La} A. Lamari, {\em Courants k\"ahl\'eriens et surfaces compactes}, Ann. Inst. Fourier (Grenoble) {\bf 49} (1999), no. 1, 263--285.
\bibitem{LY} F. Ledrappier, L.-S. Young, {\em The metric entropy of diffeomorphisms. I. Characterization of measures satisfying Pesin's entropy formula}, Ann. of Math. (2) {\bf 122} (1985), no. 3, 509--539.
\bibitem{LY2} F. Ledrappier, L.-S. Young, {\em The metric entropy of diffeomorphisms. II. Relations between entropy, exponents and dimension}, Ann. of Math. (2) {\bf 122} (1985), no. 3, 540--574.
\bibitem{Maz} B. Mazur, {\em The topology of rational points}, Experiment. Math. {\bf 1} (1992), no. 1, 35--45.
\bibitem{McM} C.T. McMullen, {\em Dynamics on $K3$ surfaces: Salem numbers and Siegel disks}, J. Reine Angew. Math. {\bf 545} (2002), 201--233.
\bibitem{McM2} C.T. McMullen, {\em Algebra and Dynamics}, Course Notes, Harvard University, 2003. \url{http://people.math.harvard.edu/~ctm/home/text/class/harvard/275/03/html/home/pn/course.pdf}
\bibitem{McM3} C.T. McMullen, {\em Dynamics on blowups of the projective plane}, Publ. Math. Inst. Hautes \'Etudes Sci. {\bf 105} (2007), 49--89.
\bibitem{Mi} Y. Miyaoka, {\em K\"ahler metrics on elliptic surfaces}, Proc. Japan Acad. {\bf 50} (1974), 533--536.
\bibitem{Ne} S.E. Newhouse, {\em Continuity properties of entropy}, Ann. of Math. (2) {\bf 129} (1989), no. 2, 215--235.
\bibitem{PS} I.I. Pjatecki\u{i}-\v{S}apiro, I.R. \v{S}afarevi\v{c}, {\em Torelli's theorem for algebraic surfaces of type $K3$}, Izv. Akad. Nauk SSSR Ser. Mat. {\bf 35} (1971), 530--572.
\bibitem{SV} N. Sibony, M. Verbitsky, in preparation, \url{http://verbit.ru/MATH/TALKS/Rigid-currents-NYU-2019.pdf}
\bibitem{Si} Y.-T. Siu, {\em Every $K3$ surface is K\"ahler}, Invent. Math. {\bf 73} (1983), no. 1, 139--150.
\bibitem{ST0} J. Song, G. Tian, {\em The K\"ahler-Ricci flow on surfaces of positive Kodaira dimension}, Invent. Math. {\bf 170} (2007), no. 3, 609--653.
\bibitem{Tod} A.N. Todorov, {\em Applications of the K\"ahler-Einstein-Calabi-Yau metric to moduli of $K3$ surfaces}, Invent. Math. {\bf 61} (1980), no. 3, 251--265.
\bibitem{To0} V. Tosatti, {\em Limits of Calabi-Yau metrics when the K\"ahler class degenerates}, J. Eur. Math. Soc. (JEMS) {\bf 11} (2009), no. 4, 755--776.
\bibitem{To2} V. Tosatti, {\em Degenerations of Calabi-Yau metrics}, in {\em Geometry and Physics in Cracow,} Acta Phys. Polon. B Proc. Suppl. {\bf 4} (2011), no. 3, 495--505.
\bibitem{To3} V. Tosatti, {\em Calabi-Yau manifolds and their degenerations}, Ann. N.Y. Acad. Sci. {\bf 1260} (2012), 8--13.
\bibitem{To4} V. Tosatti, {\em Collapsing Calabi-Yau manifolds}, preprint, arXiv:2003.00673.
\bibitem{Weh} J. Wehler, {\em $K3$-surfaces with Picard number 2}, Arch. Math. (Basel) {\bf 50} (1988), no. 1, 73--82.
\bibitem{We} A. Weil, {\em Final report on contract AF 18(603)-57}, in {\em \OE uvres scientifiques/collected papers. II. 1951--1964}, 390--395, Springer Collected Works in Mathematics. Springer, Heidelberg, 2014.
\bibitem{Ya}  S.-T. Yau, {\em On the Ricci curvature of a compact K\"ahler manifold and the complex Monge-Amp\`ere equation, I}, Comm. Pure Appl. Math. {\bf 31} (1978), 339--411.
\bibitem{Yo} Y. Yomdin, {\em Volume growth and entropy}, Israel J. Math. {\bf 57} (1987), no. 3, 285--300.
\bibitem{Zd} A. Zdunik, {\em Parabolic orbifolds and the dimension of the maximal measure for rational maps}, Invent. Math. {\bf 99} (1990), no. 3, 627--649.
\end{thebibliography}
\end{document}